\documentclass[11pt]{amsart}
\usepackage[dvipsnames]{xcolor}
\usepackage[left=0.8in,top=0.98in,right=0.8in,bottom=0.98in]{geometry}
\usepackage{mathtools, amssymb, tikz, mathrsfs}

\usepackage[T1]{fontenc}
\usepackage[colorlinks=true, linkcolor=orange, urlcolor=violet, citecolor=blue]{hyperref}
\usetikzlibrary{matrix,arrows}

\input xy 
\xyoption{all} %

\usepackage{lineno}

\usepackage{color}

\newtheorem{theorem}{Theorem}[section]
\newtheorem{lemma}[theorem]{Lemma}
\newtheorem{proposition}[theorem]{Proposition}

\newtheorem{cor}[theorem]{Corollary}
\newtheorem{conjecture}[theorem]{Conjecture}

\newtheorem{claim*}{Claim}
\newtheorem{theo}[theorem]{Theorem}
\newtheorem{definition}[theorem]{Definition}

\newcounter{cor}
\newtheorem{thmy}{Theorem}

\newcounter{theo}
\newtheorem{thmyy}{Theorem}

\theoremstyle{definition}
\newtheorem{remark}[theorem]{Remark}

\newtheorem{rmk}[theorem]{Remark}
\newtheorem{example}[theorem]{Example}

\newcommand{\G}{{\mathbb G}}

\newcommand{\C}{{\mathbb C}}

\newcommand{\Q}{{\mathbb Q}}
\newcommand{\R}{{\mathbb R}}
\newcommand{\NN}{{\mathbb N}}
\newcommand{\Z}{{\mathbb Z}}

\newcommand{\Qbar}{{\overline{\Q}}}

\newcommand{\calA}{{\mathcal A}}

\newcommand{\calC}{{\mathcal C}}

\newcommand{\calE}{{\mathcal E}}

\newcommand{\calG}{{\mathcal G}}

\newcommand{\calI}{{\mathcal I}}

\newcommand{\calO}{{\mathcal O}}

\newcommand{\calS}{{\mathcal S}}

\DeclareMathOperator{\supp}{supp}

\DeclareMathOperator{\End}{End}

\DeclareMathOperator{\ord}{ord}

\DeclareMathOperator{\Div}{Div}

\DeclareMathOperator{\Spec}{Spec}

\DeclareMathOperator{\N}{\mathbb{N}}

\DeclareMathOperator{\id}{id}

\DeclareMathOperator{\GCD}{gcd}

\numberwithin{equation}{section}
\numberwithin{table}{section}

\title{Greatest Common Divisor results on semiabelian varieties and a Conjecture of Silverman}

\author{Fabrizio Barroero}

\address{Dipartimento di Matematica e Fisica, Università degli Studi Roma 3, L.go S. L. Murialdo 1, 00146 Roma, Italy}
\email{fbarroero@gmail.com}

\author{Laura Capuano}
\email{laura.capuano@uniroma3.it}

\author{Amos Turchet}

\email{amos.turchet@uniroma3.it}

\date{}

\begin{document}

 \begin{abstract}
	A divisibility sequence is a sequence of integers $\{d_n\}$ such that $d_m$ divides $d_n$ if $m$ divides $n$. Results of Bugeaud, Corvaja, Zannier, among others, have shown that the gcd of two divisibility sequences corresponding to subgroups of the multiplicative group grows in a controlled way. Silverman conjectured that a similar behaviour should appear in many algebraic groups. We extend results by Ghioca-Hsia-Tucker and Silverman for elliptic curves and prove an analogue of Silverman's conjecture over function fields for abelian and split semiabelian varieties and some generalizations of this result. We employ tools coming from the theory of unlikely intersections as well as properties of the so-called Betti map associated to a section of an abelian scheme. 

 \end{abstract}

\maketitle

\section{Introduction}

A sequence of positive integers $\{d_n\}$ is called a \textit{divisibility sequence} if, for every $m$ that divides $n$ one has that $d_m$ divides $d_n$. Well-known examples of divisibility sequences include some linear recurrence sequences such as the Fibonacci sequence and sequences of the form $d_n=a^n-1$. In \cite[Theorem 1]{BCZ03} Bugeaud, Corvaja and Zannier proved that, given multiplicatively independent $a,b \in \Z$, for every $\varepsilon >0$ there exists a constant $c=c(a,b,\varepsilon)$ such that 
\begin{equation}\label{eq:BCZ}
	\log \GCD\{a^n-1,b^n-1\}\le \varepsilon n+c \quad \mbox{for all }n\ge 1.
\end{equation}

Recently Levin has vastly generalized this type of results in \cite{Levin2018}.

Ailon and Rudnick conjectured, see \cite[Conjecture A]{AilonRudnick}, that a stronger statement should hold, namely
\[ \GCD\{a^n-1, b^n-1\}=\GCD\{a-1, b-1\} \quad \mbox{for infinitely many } n\ge 1. \]  
This result was in fact proved by the same authors, as a Corollary of \cite[Theorem 1]{AilonRudnick}, for multiplicatively independent polynomials $a,b \in k[t]$, where $k$ is a field of characteristic zero.

Motivated by the previous results, Silverman generlized Ailon and Rudnick conjecture, considering a larger class of divisibility sequences. 
A divisibility sequence of the form $a^n-1$ corresponds to a rank $1$ subgroup of the multiplicative group $\G_m$. In this setting, the primes dividing $\GCD\{a^n-1, b^n-1\}$, and their multiplicities, can be encoded as a divisor on $\Spec \Z$. If we let $n$ vary in the positive integers, we obtain a sequence of divisors $D_{nP}$ on $\Spec \Z$ associated to the point $P=(a,b)\in \G_m^2$, where $nP=(a^n,b^n)$, called a \emph{geometric divisibility sequence}. The same construction can be performed in other algebraic groups over $\Q$, e.g. elliptic curves, obtaining sequences of divisors on $\Spec \Z$.
Silverman proposed the following conjecture:
\begin{conjecture}[{{\cite[Conjecture 10]{Silverman05}}}]\label{conj:Silverman}
	Let $\calG/\Z$ be a group scheme, let $\sigma_P \in \calG(\Z)$, and assume that
	\begin{enumerate}
		\item[(i)] the generic fiber $G = \calG \times_\Z \Q$ is an irreducible commutative algebraic group of dimension at least 2 with no unipotent part;
		\item[(ii)] the restriction $P \in G(\Q)$ of $\sigma_P$ to the generic fiber has the property that the subgroup $\Z \cdot P$ generated by $P$ is Zariski dense in $G$.
	\end{enumerate}
Then, the geometric divisibility sequence $\left( D_{nP}\right)_{n\geq 1}$ corresponding to $\sigma_P$ satisfies
	\[
		D_{nP} = D_P \qquad \text{ for infinitely many } n \geq 1.
	\]
\end{conjecture}

The conjecture is completely open in full generality. Silverman \cite{Silverman05} proved that Vojta's conjecture applied to blow ups of the generic fiber $G$ implies statements analogous to \cite[Theorem 1]{BCZ03} for divisibility sequences when $G=\G_m^2$, $G=E^2$, and $G= E \times \G_m$, for an elliptic curve $E$. 

\begin{remark}
	Corvaja and Zannier obtained results analogous to \cite{BCZ03} and \cite{Silverman05} in the function field case, such as \cite[Corollary 2.3]{CZConic}. In this setting, they prove an analogue of \eqref{eq:BCZ} with more explicit and uniform bounds, which in particular allows to recover the theorem of Ailon and Rudnick, and therefore proving some cases of the analogue of Conjecture \ref{conj:Silverman} in the function field case. For further details we refer to \cite[Section II.3]{Zanbook}.
\end{remark}
\medskip

In the spirit of Ailon and Rudnick's result \cite{AilonRudnick}, it is natural to study the conjecture in the function field case and consider group schemes $\calG \rightarrow \calC$, where $\calC$ is a nonsingular irreducible projective curve defined over $\Qbar$. If the generic fiber $G$ of $\calG$ is the product of two elliptic curves $E_1$ and $E_2$ defined over $\Qbar(\calC)$, this problem has been studied by Silverman \cite{Silverman_ff} in the isotrivial case and by Ghioca, Hsia and Tucker \cite{GHT} in general. In this setting, given a point $P=(P_1,P_2)\in G(\Qbar(\calC))$, the divisor $D_{nP}$ can be defined explicitly as
\begin{equation} \label{eq:GCD}
D_{nP}:=\sum_{\gamma \in \calC} \min\{\ord_{\gamma}(\sigma_{nP_1}^{*}(\calO_1)), \ord_{\gamma}(\sigma_{nP_2}^{*}(\calO_2))\}(\gamma).
\end{equation}
The group scheme $\calG = \calE_1 \times \calE_2 \rightarrow \calC$ is the product of two elliptic surfaces with generic fibers $E_1$ and $E_2$, and we denote by $\sigma_{nP_i}$ the section of $\calE_i \rightarrow \calC$ corresponding to the point $nP_i$ and by $\calO_i$, the image of the identity section of $\calE_i \rightarrow \calC$ for $i=1,2$. We can view $D_{nP}$ as the ``greatest common divisor'' of $\sigma_{nP_1}^{*}(\calO_1)$ and $\sigma_{nP_2}^{*}(\calO_2)$, which are divisors in $\Div(\calC)$; for this reason,
we will usually denote the right-hand side of \eqref{eq:GCD} as $\gcd\{\sigma_{nP_1}^{*}(\calO_1), \sigma_{nP_2}^{*}(\calO_2)\}$. 
When the subgroup generated by $P$ is Zariski-dense in $G(\Qbar(\calC))$, then the authors proved that the expected conclusion holds, i.e. $D_{nP} = D_P$ for infinitely many $n$. 

The first goal of this paper is to prove the function field case of Conjecture \ref{conj:Silverman} for split semiabelian varieties. 

\begin{theorem} \label{thm:semiAb}Let $\calG$ be a semiabelian scheme over a curve $\calC$, defined over the algebraic numbers. Let $G$ be the generic fiber of $
	\calG$ and suppose it is a split semiabelian variety of relative dimension greater than $1$. Let $P \in G(\Qbar(\calC))$ be a point in the generic fiber such that $\Z \cdot P$ is Zariski dense and let $\sigma_P: \calC \to \calG$ be the corresponding section. Then,
\begin{itemize}	
\item[(i)] there exists a divisor $D$ on $\calC$ such that, for every $n\ge 1$, one has $D_{nP}\le D$;
\item[(ii)] $D_{nP} = D_P$ for infinitely many $n\ge 1$. Moreover, the set of $n$ for which the equality holds is the complement in $\NN$ of a finite union of arithmetic progressions.
\end{itemize}
\end{theorem}
 
 In the case when $\calG$ is the product of two schemes of relative dimension 1, the definition of $D_{nP}$ is analogous to \eqref{eq:GCD}. In general, the image of the identity section is not a divisor in $\calG$ and we need to introduce a more general definition, see Definition \ref{def:divD_p}. Nonetheless the support of the divisors $D_{nP}$ are the points of $\calC$ where the section $\sigma_{nP}$ intersects the identity section.
 
In the case of the split semiabelian variety $E \times \G_m$ defined over $\Qbar(\calC)$, we answer a question of Silverman in \cite[Remark 6]{Silverman_ff}. 
\begin{cor} \label{cor:gcd_EGm}
	Let $\calE \rightarrow \calC$ be an elliptic surface defined over $\Qbar$ with generic fiber $E$, let $Q \in E(\Qbar(\calC))$ be a nontorsion point and let $\sigma_Q:\calC \rightarrow \calE$ be the corresponding section. 
	 Let $f\in \G_m(\Qbar(\calC))$, not a root of unity, and let $\sigma_f: \calC \rightarrow \G_m \times \calC$ the section given by $(f, \id)$. 
		Then the following hold:
	\begin{enumerate}
	\item[(i)] there exists an effective divisor $D \in \Div(\calC)$ such that 
	\[
		D_{nP}=\sum_{\gamma \in \mathcal C} \min \{\ord_{\gamma}(\sigma_{nQ}^*(\calO)), \ord_{\gamma}\sigma_{f^n}^{*}(\mathbf{1}_{\G_m}))\}(\gamma) \leq D
	\]
		for all positive integers $n\ge 1$, where $\calO$ denotes the image of the identity section of  $\calE \rightarrow \calC$ and $\mathbf{1}_{\G_m}$ is the image of the identity section of the constant scheme $\G_m \times \calC \rightarrow \calC$;
	\item[(ii)] $ D_{nP}=D_P$ for infinitely many $n\in \N$, and the set of $n$ for which the equality holds is the complement in $\NN$ of a finite union of arithmetic progressions. 
	\end{enumerate}
\end{cor}

We would like to point out that, also in the function field case, results of $\GCD$-type are closely related to Vojta's conjecture over function fields, see for example \cite{CZConic, CZGm, TurchetTrans, CaTur, GSW}.

\begin{remark}
We note that the conclusion (i) of Theorem \ref{thm:semiAb} is stronger than what is predicted by the analogue of Conjecture \ref{conj:Silverman}. However, this stronger statement does not hold in general for a nonsplit semiabelian variety, due to the presence of the so-called \emph{Ribet sections}. Indeed, a necessary condition for (i) of Theorem \ref{thm:semiAb} to hold is that the union over $n$ of the supports of $D_{nP}$ is a finite set. This relies on a relative Manin-Mumford statement (see Theorem \ref{thm:relmm} and Remark \ref{rem:MM}), which has been proved in a series of papers by Masser and Zannier for abelian schemes (see \cite{MasserZannier08,MasserZannier10,MZ12,MZ14a,MasserZannier15,MZ20}). In this setting, the necessary hypothesis that $\Z \cdot P$ is Zariski dense is also sufficient. For general semiabelian schemes this is not the case: Bertrand \cite{bertrand2011special} showed that there exists a section of a nonconstant extension of a CM elliptic curve by $\G_m$, which does not factor through any proper closed subgroup scheme (and therefore the above-mentioned hypothesis holds), but whose image meets the torsion points of the various fibers of the semiabelian scheme infinitely often (see also \cite{BMPZ}). 
 
 We can therefore see that, in this setting, the union of the supports of the $D_{nP}$ is infinite and thus Theorem \ref{thm:semiAb} does not hold. Note that this does not contradict the function field analogue of Conjecture \ref{conj:Silverman} because, even if we cannot control the support of $D_{nP}$ for all $n$, there might still be infinitely many values of $n$ such that $D_{nP} = D_P$.
\end{remark}

We have just mentioned that a fundamental ingredient for the proof of Theorem \ref{thm:semiAb} is a Manin-Mumford type statement. Actually, already Ailon and Rudnick's result relies on Manin-Mumford in $\G_m^2$ which was proved by Ihara, Serre and Tate \cite{Lang65}.
Its ``modular counterpart'' is a theorem of Andr\'e \cite{andre98}, which is nothing but the Andr\'e-Oort conjecture for $Y(1)^2$. Very recently Campagna and Dill \cite{CD22}, using Andr\'e's theorem, proved an analogue of Ailon and Rudnick's result (and more) where the polynomials $t^n-1$ are replaced by Hilbert class polynomials. 

\medskip

The second goal of the paper is to obtain a generalization of part (i) of Theorem \ref{thm:semiAb} for sequences associated to two distinct points $P$ and $Q$ in a group scheme as before. This amounts to replace the identity section in the definition of the divisor $D_{P}$ with the image of the section $\sigma_Q$, thus considering the divisor $D_{P-Q}$ (see Definition \ref{def:divD_p}). In other words, we study the locus where the images of $\sigma_P$ and $\sigma_Q$ intersect. This was already studied in the case $E_1\times E_2$ by Ghioca, Hsia and Tucker \cite[Theorem 1.1]{GHT}. Here we consider the general case of split semiabelian schemes.

\begin{theorem}\label{thm:PQ} Let $\calG$ be a semiabelian scheme over a curve $\calC$, defined over the algebraic numbers. Let $G$ be the generic fiber of $\calG$ and suppose it is a split semiabelian variety of dimension greater than 1. Let $P,Q \in G(\Qbar(\calC))$ be such that $nP\neq Q $ for all $n\in \Z$ and $\Z \cdot P$ is Zariski dense in $G$. Then, there exists an effective divisor $D \in \Div(\calC)$, independent of $n$, such that $D_{nP,Q} \leq D$.
\end{theorem}

Clearly, part (i) of Theorem \ref{thm:semiAb} is nothing but Theorem \ref{thm:PQ} with $Q=O$.

\begin{rmk}
If $G=\G_m^2$, the previous theorem is related to \cite[Theorem 1.3]{Ostafe}, which in this setting implies that, if $f_1,f_2,g_1,g_2$ are multiplicatively independent polynomials in $\Qbar[t]$ with no common zeros, then there exists a polynomial $h\in \Qbar[t]$ such that for every $n_1,n_2,m_1,m_2 \ge 0$ one has
\[ \GCD\{f_1^{n_1}-g_1^{m_1}, f_2^{n_2}-g_2^{m_2}\} \mid h.
\]
Our theorem in this case shows that there exists a polynomial $h\in \Qbar[t]$ such for every $n \in \NN$  
\[ \GCD\{f_1^n-g_1, f_2^n-g_2\} \mid h
\]
under the weaker hypothesis that $f^n\neq g$ for every $n \ge 0$.
\end{rmk}

The proofs of Theorem \ref{thm:PQ} and (i) of Theorem \ref{thm:semiAb} constists of two parts: first, one needs to show that the supports of all $D_{nP,Q}$ lie in a finite set (Theorem \ref{thm:support}) and then that the order at each point of the support is bounded independently of $n$ (Proposition \ref{prop:ord_boundQ}).

The first goal is achieved combining several unlikely intersection results and applying an argument already introduced in \cite{BC_ab} to deduce the main theorem of \cite{GHT} from one of the results in \cite{BC_semiab}. 

For the second, we first deal with the case $Q=O$ by generalizing an idea of Silverman \cite{Silverman_ff} in Lemma \ref{lemma:ord_bound}. We then deal with the general abelian case by linking the order at a point $\gamma$ of a divisor $D_{P}$ to the multiplicity of the Betti map relative to $P$ at $\gamma$, as done in \cite{CDMZ} for the one dimensional case. The linearity of the Betti map allows us to apply Lemma \ref{lemma:ord_bound} and obtain the desired bound independent of $n$ for the multiplicity of $D_{nP,Q}$ at any point of its support. For powers of the multiplicative group we reduce to a result of Ostafe \cite{Ostafe} and deduce the general split semiabelian case as a combination of the two bounds. 

\begin{remark}
We finally point out that, given a nonisotrivial elliptic surface $\mathcal E \rightarrow \calC$ and a point $P$ in the generic fiber $E(\Qbar(\calC))$, the order of the divisor $D_{nP}$ at a point $\gamma$ is the multiplicity of intersection of the section $\sigma_{nP}$ and the zero section in $\gamma$. This multiplicity has been studied independently by Corvaja, Demeio, Masser and Zannier \cite{CDMZ} and by Ulmer and Urz\'ua \cite{UU}, who proved that the set of $\gamma \in \calC$ such that a multiple of $\sigma_{P}$ intersects the zero section in $\gamma$ with multiplicity $\ge 2$ when $P$ is nontorsion is finite, while the set of $\gamma \in \mathcal C$ such that a multiple of $\sigma_{P}$ intersects the zero section in $\gamma$ is always infinite. This implies that the intersection between a multiple of $\sigma_P$ and the zero section is almost always transverse (for analogous results in the case of curves in $\G_m^2$ see \cite{MM}). On the other hand, for our purposes, the problem is somehow orthogonal: we need to show that, for \emph{fixed} $\gamma \in \mathcal C$, the multiplicity of the intersections of $\sigma_{nP}$ and $\sigma_{\calO}$ in $\gamma$ does not go to infinity when $n \rightarrow \infty$.
\end{remark}

    \subsection*{Acknowledgements} 
We thank Julian Demeio for discussing with us the results of \cite{CDMZ}. We thank Pietro Corvaja, Ariyan Javanpeykar, Siddarth Mathur, Joe Silverman and Umberto Zannier for comments and discussions. A.T. is partially supported by the projects PRIN2017: Advances in Moduli Theory and Birational Classification and PRIN2020: Curves, Ricci flat Varieties and their Interactions. The three authors are members of the INDAM group GNSAGA.

	\section{Setting and Notations}

	In this section we fix the setting and the notation that will be used in the rest of the paper.

	We denote by $\calC$ a nonsingular projective curve defined over $\Qbar$. For a given group scheme $\calG$ over $\calC$ defined over $\Qbar$ we will denote by $G$ its generic fiber, which is an algebraic group defined over $\Qbar(\calC)$. In this paper we will always assume that $G$ is irreducible and commutative, and we will write the operation additively, unless $G$ is a power of the multiplicative group.

	Given a point $P \in G(\Qbar(\calC))$, different from the identity of $G$, we denote by $\sigma_P: \calC \to \calG$ the corresponding section of the group scheme. Similarly, we will denote by $\calO$ the identity section of $\calG$, and, by abuse of notation, we will identify $\calO$ with its image $\calO(\calC)$.

	We are interested in studying the intersection of the image $\sigma_P(\calC)$ with the image of the identity section $\calO$. We will encode the information of this intersection in a divisor of $\calC$, that we will denote by $D_P$. Similarly, if we are given two distinct points $P,Q$ we can study the intersection of $\sigma_P(\calC)$ and $\sigma_Q(\calC)$ by means of a divisor $D_{P,Q}$.
	
	\begin{definition}
		\label{def:divD_p}
		Let $\calG,\calC,\calO$ and $\sigma_P$ as above. Then the \emph{divisor $D_P$ associated to $P$} is defined as follows: if $\sigma_{P}(\calC) \cap \calO = \emptyset$, then we set $D_P = 0$. On the other hand, if $\sigma_{P}(\calC) \cap \calO \neq \emptyset$ (note that the $\sigma_{P}(\calC)$ is always distinct from $\calO$ since $P$ is distinct from the identity), then the intersection is a proper closed subscheme of $\calC$ via the isomorphism $\sigma_P: \calC \to \sigma_P(\calC)$. Since $\calC$ is a nonsingular projective curve the subscheme is an effective divisor that we define to be $D_P$.
	
		Given another point $Q \in G(\Qbar(\calC))$, different from $P$, we define similarly $D_{P,Q}$ as $D_{P-Q}$, i.e. the divisor associated to the point $P-Q$. 
	\end{definition}

	In the case where the group scheme $\calG$ is defined over $\Z$, the previous definition was given in \cite[Definition 4]{Silverman05}. Similar to the arithmetic case of Silverman, we can make the order of the divisor $D_{P,Q} $ at a point of $\calC$ explicit.
	
	\begin{remark}\label{rmk:mult}
		Given $P,Q$ as before, with corresponding sections $\sigma_P,\sigma_Q$, let $P_\gamma = \sigma_P(\gamma)$. The support of $D_{P,Q}$ consists precisely of the $\gamma \in \calC(k)$ such that $\sigma_Q(\gamma) = \sigma_P(\gamma) = P_\gamma$. Let $\calI_P,\calI_Q$ be the ideal sheaves corresponding to $\sigma_P$ and $\sigma_Q$ in the completed local ring $\hat\calO_{\calG,P_\gamma}$. Then the order of $\gamma$ in $D_{P,Q}$ is
		\[
			\ord_\gamma D_{P,Q} = \text{length} \dfrac{\hat\calO_{\calG,P_\gamma}}{\langle \calI_P, \calI_Q \rangle}.
		\]
		In particular the divisor $D_{P,Q}$ encodes information both on the points where the two curves $\sigma_P(\calC)$ and $\sigma_Q(\calC)$ intersects, and on the multiplicity of their intersection.
	\end{remark}
	
	\begin{remark}\label{rmk:dim1} 
		When $\calG$ has relative dimension 1, Definition \ref{def:divD_p} can be made even more explicit. Indeed in this case the image of the zero section is a divisor on $\calG$. In this case the divisor $D_P$ can be defined as the pullback $\sigma_P^*(\calO)$ (and $D_{P,Q} = \sigma_P^*(\sigma_Q(\calC))$). This coincides with Definition \ref{def:divD_p} but has the disadvantage that is does not extend to higher dimension. 
	\end{remark}

	We will now make the construction of the divisor $D_P$ explicit in a couple of relevant situations, namely when $G$ is an elliptic curve or the multiplicative group. 
	
	\begin{example}
	  \label{ex:D_ell}
	Let $E$ be an elliptic curve defined over $\Qbar(\calC)$, and let $P \in E(\Qbar(\calC))$ be a point. We fix an elliptic scheme $\calE \to \calC$ with zero section $\calO$ which is a proper model of $E$ over $\calC$, and a positive integer $n$. Then we can use Remark \ref{rmk:dim1} with $\calG = \calE$ to define the divisor $D_{[n]P} $ as follows: one considers the section $\sigma_{[n]P}$ associated to the point $[n]P \in E(\Qbar(\calC))$. The section $\sigma_{[n]P}$ defines a pullback map on divisors $\sigma_{[n]P}^*: \Div (\calE) \to \Div \calC$ so that $\sigma_{[n]P}^*(\calO)$, i.e. the pullback of the zero section of the elliptic scheme $\calE \to \calC$, is a divisor on $\calC$, which we denote by $D_{[n]P}$. 
	
	If $\gamma$ is a point of $\calC$ over which $\calE$ has good reduction, then $\gamma$ is in the support of $D_{[n]P}$ if and only if $\sigma_{[n]P}(\gamma) = \calO(\gamma)$ in the elliptic curve $\calE_\gamma$. 
	\end{example}

	\begin{example}
	  \label{ex:D_f}
	  Similarly, if $f \in \G_m(\Qbar(\calC)) = \Qbar(\calC)^\times$ is a rational function, we can use Remark \ref{rmk:dim1} to define the divisor $D_{f^n} = \sigma_{f^n}^*(\mathbf{1}_{\G_m})$ for every positive integer $n$. Here $\sigma_{f^n}$ denotes the section of $\calG$ corresponding to the rational function $f^n$ and $\mathbf{1}_{\G_m}$ is the identity section corresponding to the point $1 \in \G_m(\Qbar(\calC))$. The construction is completely analogous to the case of Example \ref{ex:D_ell}. 
	
	In this setting, a point $\gamma \in \calC(\Qbar)$ is in the support of $D_{f^n}$ if and only if $\sigma_{f^n}(\gamma) = 1$ in $\G_m = \G_m \times \{ \gamma \}$, i.e. if $f^n(\gamma) = 1$.
	\end{example}

	We conclude this section by stressing the link between the divisor $D_{nP}$ and the $\GCD$ problems mentioned in the introduction.
	
	\begin{example}\label{ex:D_gcd}
		Let $G$ be the semiabelian variety $E \times \G_m$, where $E$ is some elliptic curve defined over $\Qbar(\calC)$. Then, the divisor $D_P$ associated to a point of $G$ can be expressed as a $\GCD$ of the divisors defined in Examples \ref{ex:D_ell} and \ref{ex:D_f}. Indeed, given a nonzero point $Q \in E(\Qbar(\calC))$, a function $f \in \G_m(\Qbar(\calC))\setminus \{1\}$, and the point $P = (Q,f) \in G(\Qbar(\calC))$, Definition \ref{def:divD_p} yealds three divisors $D_Q, D_f$ and $D_P$ in $\calC$. Then, it is easy to check that $D_P = \GCD(D_Q, D_f)=\sum_{\gamma \in \mathcal C} \min \{\ord_{\gamma}(\sigma_{Q}^*(\calO)), \ord_{\gamma}\sigma_{f}^{*}(\mathbf{1}_{\G_m}))(\gamma) $ (cfr \cite[Proposition]{Silverman05} for an analogue over the integers).
	\end{example}

\section{Unlikely Intersections and finiteness of the support} 

We start this section by formulating two theorems that are going to be the key ingredient to control the support of the divisors appearing in Theorems \ref{thm:semiAb} and \ref{thm:PQ}. The first is a consequence of Pink's conjecture (Conjecture 6.1 in \cite{PinkUnpubl}) and is obtained as a combination of several results.
 
\begin{theorem}\label{th:unl1}
	Let $\calC$ be a nonsingular projective curve defined over $\Qbar$, let $\calG \rightarrow \calC$ a group scheme defined over $\Qbar$ and denote by $G$ its generic fiber. Suppose that $G$ is a split semiabelian variety of dimension at least 2. Let $P\in G(\Qbar(\calC))$ and $\sigma_{P}$ be the corresponding section. If the set
		\begin{equation}\label{Eq:set}
		 \{ \gamma \in \calC(\Qbar) : \text{there exists an algebraic subgroup $H$ of $G$ of codimension $\geq 2$ such that } \sigma_P(\gamma) \in H_\gamma \}
	\end{equation}
is infinite, then $P$ lies in a proper algebraic subgroup of $G$, i.e., $\Z\cdot P$ is not Zariski-dense in $G$. 
\end{theorem}

\begin{proof}
	For any finite cover $\calC'\to \calC$, we have a semiabelian scheme $\mathcal{G}'=\mathcal{G}\times_{\calC} \calC'$ over $\calC'$ whose generic fiber is the base-change $G_{\Qbar(\calC')}$ of $G$. The section $\sigma_{P}$ extends to  $\sigma'_{P}:\calC' \to \mathcal{G}'$. Now, the set \eqref{Eq:set} is infinite if and only if the same set with $\mathcal{C}'$ and $\sigma'_{P}$ in place of $\mathcal{C}$ and $\sigma_{P}$ is infinite. Therefore, we are allowed to assume that all semiabelian subvarieties and endomorphisms of $G$ are defined over $K$.
	
	As our statement is invariant under isogenies, we may moreover assume that $G= \prod_{i=1}^{d}A_i^{e_i}$ for nonnegative integers $d,e_1, \dots, e_d$ where the $A_i$ are pairwise nonisogenous simple abelian varieties or $\G_m$. As we are considering algebraic subgroups of $G$ of codimension $\geq 2$, we may restrict ourselves to the cases and corresponding results listed below and obtain our claim.
	\begin{enumerate}
		\item $G=\G_m^n$: \cite{Maurin} (see also \cite{BMZ99,BHMZ,CMPZ}).
		\item $G=A$ an isoconstant abelian variety: \cite{HP} after previous partial results  \cite{RemVia,Ratazzi08,Carrizosa,Viada2008,galateau2010}.
		\item $G=A$ a nonisoconstant abelian variety: \cite{BC_ab} and the earlier \cite{BC_ell,BC_semiab}.
		\item $G$ an isoconstant semiabelian variety: \cite{BKS} (see also \cite{BarSha,barroero2021multiplicative}).
		\item $G=E^n\times \G_m^l$ for a nonisoconstant elliptic curve $E$: \cite{BC_semiab}.\qedhere
	\end{enumerate}
\end{proof}

The following theorem is usually considered as a relative version of the Manin-Mumford conjecture. Indeed, if $G$ is isoconstant, then this is implied by the usual Manin-Mumford, proved by Laurent \cite{Laurent84}, Raynaud \cite{Raynaud} and Hindry \cite{Hindry1988} for $\G_m^n$, abelian varieties and semiabelian varieties respectively. Recently Masser and Zannier have investigated the relative case for abelian varieties in a series of papers \cite{MasserZannier08,MasserZannier10,MZ12,MZ14a,MasserZannier15,MZ20}. All these works imply the following theorem, which is also implied by Theorem \ref{th:unl1}. We state it separately because in same cases one does not need to use the full strength of the above theorem but a relative Manin-Mumford statement is enough (see Remark \ref{rem:MM}).

\begin{theorem}\label{thm:relmm}
	Let $\calC$ be a nonsingular projective curve defined over $\Qbar$, let $\calG \rightarrow \calC$ a group scheme defined over $\Qbar$ and denote by $G$ its generic fiber.
	Suppose that $G$ be a split semiabelian variety over $\Qbar(\calC)$ of dimension at least 2, and let $P\in G(\Qbar(\calC))$ and $\sigma_{P}$ be the corresponding section. If there are infinitely many $\gamma\in \calC(\Qbar)$ such that $\sigma_{P}(\gamma)$ is torsion in $G_\gamma$, then $P$ lies in a proper algebraic subgroup of $G$, in particular $\Z\cdot P$ is not Zariski-dense in $G$. 
\end{theorem}

We will now apply Theorems \ref{th:unl1} and \ref{thm:relmm} to prove that the following result will allow us to control the support of the divisors we are studying.

\begin{theorem}\label{thm:support}
	Let $\calC$ be a nonsingular projective curve defined over $\Qbar$, let $\calG \rightarrow \calC$ a group scheme defined over $\Qbar$ and denote by $G$ its generic fiber. Suppose that $G$ is a split semiabelian variety over $\Qbar(\calC)$, and let $P,Q \in G(\Qbar(\calC))$. 
If the set
	\begin{equation}\label{Eq:bigunion}
	\bigcup_{n \geq 1} \{ \gamma \in \calC(\Qbar) : \sigma_{nP}(\gamma) = 
	\sigma_{Q}(\gamma) \}
\end{equation}
 is infinite then at least one of the following holds:
 \begin{enumerate}
 	\item $nP=Q$ for some integer $n$, or
 	\item there exists a finite cover $\calC'\to \calC$ and an isogeny $\alpha : G_{\Qbar(\calC')} \to G_1\times G_2$ such that $\Z\cdot \pi_1(P)$ is Zariski dense in $G_1$, $	\pi_2(P)=\pi_2(Q)=O_{G_2}$ and $\dim (G_1)=1$, where $\pi_i$ is the composition of $\alpha$ with the projection on $G_i$. 
 \end{enumerate}
\end{theorem}
\begin{proof} We assume \eqref{Eq:bigunion} is infinite and (1) is false. We want to prove that (2) holds.
	 
	For any finite cover $\calC'\to \calC$ we have a semiabelian scheme $\mathcal{G}'=\mathcal{G}\times_{\calC} \calC'$ over $\calC'$ whose generic fiber is the base-change $G_{\Qbar(\calC')}$ of $G$. The sections $\sigma_{nP}$ and $\sigma_Q$ extend to sections $\sigma'_{nP} ,\sigma'_Q:\calC' \to \mathcal{G}'$. Now, clearly \eqref{Eq:bigunion} is finite if and only if $\bigcup_{n \geq 1} \{ \gamma \in \calC'(\Qbar) : \sigma_{nP}'(\gamma) = \sigma_{Q}'(\gamma) \} $ is finite. We may and will replace $\calC$ by $\calC'$ tacitly and assume that the morphisms of algebraic groups we are considering are defined over $\Qbar(\calC)$.
	
	By considering the Zariski-closure of $\Z\cdot P$ we may always find an isogeny $\alpha : G\to G_1\times G_2$ such that $\Z\cdot \pi_1(P)$ is Zariski dense in $G_1$ and $\pi_2(P)=O_{G_2}$.

	Now, possibly extending $\Qbar(\calC)$ again we have an isogeny $\beta:G_1\to H_1\times \dots \times H_r$ where the $H_j$ are geometrically simple factors (geometrically simple abelian varieties or $\G_m$). As before, we may assume that they are all defined over $\Qbar(\calC)$.

We may also spread out and obtain semiabelian schemes $\mathcal{G}_1, \mathcal{G}_2,\mathcal{H}_1,\dots , \mathcal{H}_r$ over $\calC$ whose generic fibers are $G_1,G_2, H_1, \dots , H_r$ respectively, and morphisms $\mathcal{G}\to\mathcal{G}_1\times_{\calC} \mathcal{G}_2$ and $\mathcal{G}_1\to\mathcal{H}_1\times_{\calC}\dots \times_{\calC} \mathcal{H}_r$ over $\calC$.

We let $G_\gamma$ be the fiber of $\mathcal{G}$ above $\gamma \in \calC$, and we use the same notation for the $ G_i$ and the $H_j$ and for the morphisms between them.	

Now, since $	\pi_2(P)=O_{G_2}$ and \eqref{Eq:bigunion} is infinite, we must have that $\pi_2(Q)=O_{G_2}$.
Also, we clearly have that, if $\bigcup_{n \geq 1} \{ \gamma \in \calC(\Qbar) :(\pi_1)_\gamma( \sigma_{nP}(\gamma)) =(\pi_1)_\gamma( \sigma_{Q}(\gamma) )\}$ is finite, then \eqref{Eq:bigunion} is finite. We may then assume $G=G_1$ and $\Z\cdot P$ is Zariski-dense in $G$.

Now, we note that we may assume $G= H_1\times \dots \times H_r$. Indeed, the isogeny $\beta$ preserves the 2 conditions above and \eqref{Eq:bigunion} is infinite if and only if $\bigcup_{n \geq 1} \{ \gamma \in \calC(\Qbar) :(\beta)_\gamma( \sigma_{nP}(\gamma)) =(\beta)_\gamma( \sigma_{Q}(\gamma) )\}$ is infinite.

Since we are not in case (1), there exists a sequence $(\gamma_m)_{m\in \N}$ of distinct points of $ \calC(\Qbar)$ and a sequence of integers $(n_m)_{m\in \N}$ with $|n_m|\rightarrow \infty$ such that $ \sigma_{n_mP}(\gamma_m) = 
\sigma_{Q}(\gamma_m)$. 

If we are not in case (2), then we have a factor $H$ of $G$, that is either
\begin{enumerate}
	\item[(i)] geometrically simple with $\dim(H)\geq 2$, or
	\item[(ii)] a product of two one dimensional semiabelian varieties.
\end{enumerate} 

Let $p:G\to H$ be the projection. If we show that there is an $n\in \N$ with $np(P)=p(Q)$ then we are done. Indeed, we have two relations 
$$
p_{\gamma_m}(\sigma_{n_mP}(\gamma_m) )=p_{\gamma_m}( 
\sigma_{Q}(\gamma_m) )\ \ \text{ and }\ \ p_{\gamma_m}(\sigma_{nP}(\gamma_m) )=p_{\gamma_m}( 
\sigma_{Q}(\gamma_m))
$$
which, as $|n_m|\rightarrow \infty$, are independent for $m$ large enough. This implies that $p_{\gamma_m}(\sigma_{P}(\gamma_m) )$ is torsion for infinitely many $\gamma_m$. By Theorem \ref{thm:relmm}, this contradicts the fact that $\Z\cdot P$ is Zariski-dense in $G$.

If we are in case (i), then we can just apply Theorem 8.5 of \cite{BC_ab} to get our sought relation. We can then assume that we are in case (ii), hence we can write $H=H_1\times H_2$ and let $p_j:G\to H_j$ be the projection. We may also assume, arguing as above, that any isogeny between $H_1$ and $H_2$ is defined over $\Qbar(\calC)$, as well as any endomorphism of $H_1$ and $H_2$.

Now, to conclude, one could simply imitate the proof of Theorem A-1 of \cite{BC_ab} and invoke the right unlikely intersection result. We include the proof anyway for the reader's convenience. 

We let $P_j=p_j(P)$ and $Q_j=p_j(Q)$. We will moreover write $P_j(\gamma)$ for $p_j(\sigma_{P}(\gamma))$, and similarly $Q_j(\gamma)=p_j(\sigma_{Q}(\gamma))$.

We have
 \begin{equation}\label{eqa1}
	{n_{m} }P_1(\gamma_m)=Q_1(\gamma_m),
\end{equation} and \begin{equation}\label{eqa11}
{n_{m}}	P_2(\gamma_m)=Q_2(\gamma_m).
\end{equation} 

Suppose first that $H_1$ and $H_2$ are not isogenous. Then, by Theorem \ref{th:unl1}, $(P_1,Q_1,P_2,Q_2)$ lies in a proper algebraic subgroup of $H_1^2\times H_2^2$. We may assume by symmetry that there are $a_1,b_1\in \End(H_1)$ such that $a_1P_1=b_1Q_1$. Note that $b_1=0 $ would contradict the fact that $\Z\cdot P$ is Zariski-dense in $G$. This relation, together with \eqref{eqa1}, implies that $P_1(\gamma_m)$ is torsion for $m$ large enough. Since $P_1$ cannot be identically torsion, Theorem \ref{th:unl1} gives us $a_2P_2=b_2Q_2$ for some $a_2,b_2\in \End(H_2)$ with $b_2\neq 0$ as before. If we combine this with \eqref{eqa11} we have that  $P_2(\gamma_m)$ is torsion for $m$ large enough. We get a contradiction by applying Theorem \ref{thm:relmm}.

We are then left to deal with the case in which $H_1$ and $H_2$ are isogenous. This is similar to what we have just done in the nonisogenous case but a bit more involved. We may assume $H_1=H_2=:H$.

This time Theorem \ref{th:unl1} gives $a_1,a_2,b_1,b_2\in \End(H)$, not all zero such that 
\begin{equation}\label{eqa4}
	a_1P_1+a_2P_2=b_1Q_1+b_2Q_2.
\end{equation}
As before, $b_1,b_2$ cannot be both zero since otherwise this would contradict the fact that $\Z\cdot P$ is Zariski-dense in $G$. Therefore, we may assume $b_1\neq 0$. Then, combining \eqref{eqa1} and \eqref{eqa4} we have that $P_1(\gamma_m),P_2(\gamma_m)$ and $Q_2(\gamma_m)$ satisfy nontrivial dependence relations for $m$ large enough. As the coefficient in front of $P_2(\gamma_m)$ is independent of $\gamma_m$, this relation is independent of \eqref{eqa11}. Thus, by Theorem \ref{th:unl1}, we have
\begin{equation}\label{eqa5}
	{c_1} P_1 +{c_2}P_2=d_2Q_2,
\end{equation}
for $c_1, c_2, d_2\in \End(H)$, not all zero, and we may assume $d_2\neq 0$ as before. Combining \eqref{eqa11} and \eqref{eqa5}, we have
\begin{equation}\label{eqa6}
	c_1P_1(\gamma_m) = ( d_2 n_{m}-c_2 )P_2(\gamma_m).
\end{equation}
This, together with \eqref{eqa1}, gives
\begin{equation*}
	{e_1} P_1 +{e_2}P_2=f_1Q_1,
\end{equation*}
for some $e_1,e_2, f_1\in \End(H)$ with $f_1\neq 0$.
Combining this with \eqref{eqa1} again we obtain 
\begin{equation}\label{eqa8}
{(e_1 -f_1 n_{m} )}	P_1(\gamma_m) = e_2 P_2(\gamma_m).
\end{equation}
Now, \eqref{eqa6} and \eqref{eqa8} cannot be independent relations because of Theorem \ref{th:unl1}. Therefore,
$$
c_1 e_2 =(e_1 -f_1 n_{m} )( d_2 n_{m}-c_2 )
$$
that, since $f_1d_2\neq0$, gives a nontrivial quadratic equation for $n_{m}$ independent of $m$. This is impossible because $|n_m|\rightarrow \infty$.
	\end{proof}

\begin{remark}\label{rem:MM}
	One can see that, in the case $Q = O$ (or more generally $Q$ is torsion) there is no need to invoke Theorem \ref{th:unl1} but only Theorem \ref{thm:relmm} in the proof of Theorem \ref{thm:support}. \end{remark}

	\section{Bounding the multiplicity}\label{sec:mult}
In this section we prove that the order of the divisors $D_{nP}$ and $D_{nP,Q}$ at points of their support is bounded independently of $n$.

The first case generalizes an argument given by Silverman in the case of elliptic curves.

    \begin{lemma}[cfr.~{{\cite[Lemma 4 and Remark 2]{Silverman_ff}}}]
	    Let $G$ be a connected algebraic group defined over the function field $\Qbar(\calC)$ of a nonsingular projective curve $\calC$, let $\calG$ be a model of $G$ over $\calC$ and let $\calO$ be the image of the section of $\calG$ corresponding to the identity element of $G$. Let $P$ be a nontorsion point of $G$ and let $D_P$ be the divisor on $\calC$  defined according to Definition \ref{def:divD_p}. Let $\gamma$ be a point of $\calC$. Then,
	    \begin{enumerate}
		    \item if $\ord_\gamma D_P \geq 1$ 
 then
			    \begin{equation*}
				    \ord_\gamma D_{nP} = \ord_\gamma D_P \qquad \text{ for all } n \neq 0;
			    \end{equation*}
		    \item there is an integer $m = m(\calG,P,\gamma)$ so that
			    \begin{equation*}
	\ord_\gamma D_{nP} \in \{ 0, m \} \qquad \text{ for all } n \neq 0,
\end{equation*}
in particular $\ord_\gamma D_{nP}$ is bounded independently of $n$.
	    \end{enumerate}
	    \label{lemma:ord_bound}
    \end{lemma}
    \begin{proof}
	    The proof is analogous to \cite[Lemma 4]{Silverman_ff}. If we denote by $\mu_r: \calG \to \calG$ the multiplication by $r$ map, then the subscheme $\mu_r^*\calO$ is the union of $\calO$ and the set $T_r$ of nonzero $r$-torsion sections of $\calG$. Since $G$ is a connected algebraic group in characteristic zero, $\calG$ is divisible which implies that $\mu_r$ is \'etale in a neighborhood of the identity element $\calO_\gamma$ of $\calG_\gamma$ hence $\calO \cap T_r = \emptyset$. This implies that
	    \begin{equation*}
		    \ord_\gamma D_{nP} = \ord_\gamma D_P,
	    \end{equation*}
	    proving (1).

	    To prove (2) we can assume that there exist at least two multiples $n_1P$ and $n_2P$ such that $\ord_\gamma D_{n_i P} \geq 1$ for $i=1,2$, otherwise there is nothing to prove. Then (1) implies 
	    \[
		    \ord_\gamma D_{n_1P} = \ord_\gamma D_{n_1 n_2 P} = \ord_\gamma D_{n_2P}=: m.
	    \]
	    This implies that for every $n \neq 0$ either $\ord_\gamma D_{nP} = 0$ or $\ord_\gamma D_{nP} = \ord_\gamma D_{n_1P} = m$ as wanted.
	\end{proof}
	
	The second statement we obtain is the analogous bound for divisors of the form $D_{nP,Q}$. This was obtained in \cite{GHT} for products of elliptic curves; however, we use a completely different argument to deduce the bound for split semiabelian varieties of the form $A \times \G_m^\ell$.

	\begin{proposition}
		\label{prop:ord_boundQ}
		Let $A$ be an abelian variety over $\Qbar(\calC)$ and let $G = A \times \G_m^\ell$ for some $\ell \geq 0$. Let $P,Q$ be two points of $G$ such that $Q \neq nP$ for every $n \geq 1$. For every $n\ge 1$ let $D_{nP,Q}$ be the divisor on $\calC$  defined according to Definition \ref{def:divD_p}. Then, there exists an integer $m = m(G,P,Q,\gamma)$ such that $\ord_\gamma D_{nP,Q} \leq m$ for every $n\ge 1$.
	\end{proposition}
	
	In order to prove Proposition \ref{prop:ord_boundQ}, we first deal with the case of abelian schemes. In this case we can characterize the order of the divisor at a point $\gamma$ in terms of the Betti map, following \cite[Lemma 2.7]{CDMZ}. 

	The Betti map has been a central tool in the study of Diophantine problems, and recently its rank has been investigated by several authors, see for example \cite{ACZ,Gao,DGH} and references therein. We will begin by recalling basic facts about the Betti map.
\medskip

Given a point $P$ in a complex abelian variety $A$ of dimension $g$, its abelian logarithm can be expressed as a linear combination of the periods of $A$ with real coefficients, usually called the \emph{Betti coordinates} of $P$. In the case of an abelian scheme $\mathcal A \rightarrow \mathcal C$ and a section $\sigma_P:\mathcal C \rightarrow \mathcal A$, these coordinates become a system of multivalued real-analytic functions. %

More precisely, since $\calC$ is smooth, by a result of Raynaud the abelian scheme $\pi: \calA \to \calC$ carries a principal polarization and has a level $\ell \geq 3$ structure. By \cite[Proposition 2.1]{DGH}, for every $c_0 \in \calC^{an}(\C)$ there exists a connected open neighborhood $\Delta$ of $c_0$ and a real-analytic map $b_\Delta: \calA_\Delta := \pi^{-1}(\Delta) \to (\R/\Z)^{2g}$ which is fiberwise a group isomorphism.

Explicitly, if $\Delta$ is simply connected, for each $c \in \Delta$ one can define a basis of the period lattice at each fiber, $\rho_1(c),\dots,\rho_{2g}(c)$, as holomorphic functions of $c$. In the identification $\calA_c = \pi^{-1}(c)$ with $\C^g / \Lambda_c$, with $\Lambda_c$ the lattice generated by the periods at $c$, each $x \in \calA_c(\C)$ is the class of
\[
	b_1(x)\rho_1(c) + \dots + b_{2g}(x)\rho_{2g}(c)
\]
for some real numbers $b_1(x),\dots,b_{2g}(x)$. Then $b_\Delta(x)$ is in the class of the $2g$-tuple $(b_1(x),\dots,b_{2g}(x))$ modulo $\Z^{2g}$.
If $\sigma_P: \calC \to \calA$ is a section,  the composition $\beta_P:= \beta_\Delta \circ \sigma_P \vert_{\Delta}$, given in coordinates by 
 \begin{align*}
	\beta_P: \Delta &\longrightarrow (\R/\Z)^{2g}	\\
	c &\longmapsto \left(\beta_{1,P}(c),\dots,\beta_{2g,P}(c)\right),
 \end{align*} 
is called the \emph{Betti map associated to $\sigma_P$}, with respect to the neighborhood $\Delta$. 

\medskip
We denote by $m_{\beta_P}(\gamma)$ the multiplicity of $\beta_P$ at the point $\gamma$. In the following lemma we show that $m_{\beta_P}(\gamma)$ coincides with $\ord_\gamma D_P$ in an abelian scheme.

	\begin{lemma}
		\label{lemma:Betti_ab}
		Let $\sigma_P: \calC \to \calA$ be a nontorsion section of an abelian scheme $\calA$ and let $\gamma \in \calC$ such that $\sigma_P(\gamma) = \calO(\gamma)$ in $\calA_\gamma$. Then the order $\ord_\gamma D_P$ equals $m_{\beta_P}(\gamma)$.
	\end{lemma}

	\begin{proof}
		The statement is local, so we can assume that, in a small neighborhood of $\gamma = 0$, the sections $\sigma_P$ and the zero section $\calO$ are locally given respectively by
		\[
			\sigma_P: t \mapsto (f_1(t),\dots,f_g(t),t) \qquad \calO: t \mapsto (0,\dots,0,t), 
		\]  
		where $f_j(t)$ are complex analytic functions for every $j=1, \ldots, g$.
		In this setting, the order $m$ of $D_P$ at $0$ is given by
		\[
			m := \ord_0 D_P = \min \left\{i: \dfrac{d^{(i)} f_j}{d t}(0) \neq 0 \text{ for some } j = 1,\dots,g \right\}.
		\]
		Similarly the multiplicity of the Betti map is given by
		\begin{equation}\label{eq:mbeta}
			m_{\beta_P} = m_{\beta_P}(0) = \min\left\{i: \dfrac{d^{(i)}\beta_{j,P}}{d x^{(i)}}(0) \neq 0, \text{ for some } j =1,\dots,g \right\}.
		\end{equation}

		We note that, in our setting, both minima are strictly positive since $\sigma_P(0) = \calO(0)=(0,\dots,0)$.
		
		Let $\tilde{\sigma}_P$ be an abelian logarithm of $\sigma_P$, which locally in a neighborhood $U$ of $\gamma = 0$ is given by
		\begin{equation} \label{eq:sigma_int}
			\tilde{\sigma}_P(u) = \left( \int_{\calO(u) \to \sigma_P(u)} \omega_1, \dots, \int_{\calO(u) \to \sigma_P(u)} \omega_g \right),
		\end{equation}
		where $\omega_1,\dots,\omega_g$ are a basis of $\Omega^1_{\calA/\calC}$ and we are integrating on some choice of path from $\calO(u)$ to $\sigma_P(u)$. In our setting, the sheaf of one forms of $\calA$ is locally generated by $d
		x_1,\dots,dx_g$, where the $x_1, \ldots, x_{g}$ are local parameters. Therefore we can express locally each $\omega_i$ as the sum 
		\[
			 \sum_{j=1}^g s_{i,j}(x_1,\dots,x_g) dx_j,
		\]
		where the $s_{i,j}$ are power series in $x_1, \ldots, x_g$ and the determinant of the $g \times g$ matrix $S(x_1, \ldots, x_g):=(s_{i,j}(x_1, \ldots, x_g))$ does not vanish when evaluated in $\underline{0}$.
		 Notice that a path from $\calO(u)$ to $\sigma_P(u)$ on $\calA_u$ corresponds, via the local parameters, to a path in $\C^g$ from $\underline{0}$ to $\underline{f}(u):=(f_1(u), \ldots, f_g(u))$. 
		Then, a similar argument as in \cite[Lemma 2.7]{CDMZ} 
		allows to prove that $\Vert \tilde{\sigma}_P(u) \Vert \sim \Vert u \Vert^m$. Indeed, computing the integrals defining $\tilde{\sigma}_P(u)$ in \eqref{eq:sigma_int} using the power series expansions of the $s_{ij}$, one has that
        \[
			\Vert \tilde{\sigma}_P(u) \Vert= \left \Vert S(\underline 0)\cdot \underline{f}(u)^t \right \Vert +O(\Vert \underline{f}(u) \Vert^2),
	     \]
		 which implies that 
		 \[
			\Vert \underline{f}(u) \Vert \ll \Vert \tilde{\sigma}_P(u) \Vert \ll \Vert \underline{f}(u) \Vert.
		\]
		Moreover, since the $f_j(t)$ are complex analytic functions and $u$ is in a neighborhood of $0$, one has that $\Vert \underline{f}(u) \Vert \sim \Vert u \Vert^m$, giving $\Vert \tilde{\sigma}_P(u) \Vert \sim \Vert u \Vert^m$ as wanted.

		On the other hand, by the same argument as in the end of \cite[Lemma 2.7]{CDMZ}, and the definition of the Betti map, we have that
		\[
			\Vert \tilde{\sigma}_P(u) \Vert \sim \Vert\left(\beta_{1,P}(u), \dots, \beta_{2g,P}(u)\right)\Vert.
		\]
		Since the $\beta_i$ are real analytic we have that
		\[
			\Vert \left(\beta_{1,P}(u), \dots, \beta_{2g,P}(u)\right)\Vert \sim \Vert u \Vert^{m_{\beta_P}}.
		\]
		Combining the last two estimates with the fact that $\Vert \tilde{\sigma}_P(u) \Vert \sim \Vert u \Vert^m$, we obtain $m=m_{\beta_P}$, as wanted.
	\end{proof}

	Using the relation with the Betti map we can give a uniform bound for the order of a divisor of the form $D_{nP,Q}$ at a point $\gamma$ in its support in the case of abelian schemes.

	\begin{lemma}\label{lemma:mult_ab}
		Let $\sigma_P,\sigma_Q$ be two sections of an abelian scheme $\calA \to \calC$ defined over $\Qbar$ corresponding to points $P,Q \in A(\Qbar(\calC))$, such that $nP \neq Q$ for every $n\in \Z$, and let $\gamma \in \calC$. Then, there exists $m=m(\calA,P,Q,\gamma)$ such that $\ord_\gamma D_{nP,Q} \leq m$ for every $n \geq 1$.
	\end{lemma}
	\begin{proof}

	We can assume $Q \neq O$ since otherwise the conclusion follows by Lemma \ref{lemma:ord_bound}. 
	Furthermore, if there is at most one $n$ such that $\ord_{\gamma}D_{nP,Q}\ge 1$, then there is nothing to prove, so we will assume there exists more than one $n$ satisfying this inequality. 
	This implies that $\sigma_Q(\gamma)$ is a torsion point and therefore there exists a positive integer $a$ such that $\sigma_{aQ}(\gamma) = \sigma_{anP} (\gamma) = \mathcal{O}(\gamma)$. Note that, by Lemma \ref{lemma:ord_bound} and Lemma \ref{lemma:Betti_ab}, we have that $m_{\beta_{a(nP-Q)}}(\gamma)  = m_{\beta_{nP-Q}}(\gamma) $.

		From the local definition of the multiplicity of the Betti map in \eqref{eq:mbeta} and the fact that
		\[
		\beta_{P_1+P_2} = \beta_{P_1} + \beta_{P_2},
		\]
		we get that	
		\begin{equation}\label{eq:betti1}
			m_{\beta_{nP - Q}}(\gamma) = m_{\beta_{a(nP-Q)}}(\gamma) \geq \min \{ m_{\beta_{anP}}(\gamma) , m_{\beta_{aQ}}(\gamma) \}.
		\end{equation}
		On the other hand, there exists at most one $n$ such that \eqref{eq:betti1} is not an equality (namely, when the value of $n$ makes the linear combination of the corresponding derivatives vanish at $\gamma$). This implies that there exists a constant $m'=m'(\calA,P,Q,\gamma)$ {independent of $n$} such that
		\begin{equation}
			\label{eq:betti2}
			m_{\beta_{nP-Q}}(\gamma) \leq \min \{ m_{\beta_{anP}}(\gamma) , m_{\beta_{aQ}}(\gamma)\} + m'.
		\end{equation}
		Now, by Lemma \ref{lemma:Betti_ab}, we know that the multiplicities of the Betti map equal the orders of the corresponding divisors, thus we obtain from \eqref{eq:betti2} that
		\begin{equation*}
			\ord_\gamma D_{nP,Q} \leq \min \{ \ord_\gamma D_{anP} , \ord_\gamma D_{aQ}\} + m'.
		\end{equation*}
		To conclude, we apply Lemma \ref{lemma:ord_bound} to bound $\ord_\gamma D_{anP}$ independently of $n$ and we note that both $\ord_\gamma D_{aQ}$ and $m'$ do not depend on $n$, thus finishing the proof.
	\end{proof}

	The same strategy would allow to show that the order of $D_{nP-Q}$ can be characterized in terms of Betti map also in the multiplicative group. In this case, however, we can obtain directly that the order is bounded independently of the integer $n$.

	\begin{lemma}\label{lemma:mult_Gm}
	Let $P,Q$ be two points in $\G_m^\ell$ such that $Q$ is nontorsion and $Q \neq nP$ for every $n \geq 1$. Let $\sigma_P,\sigma_Q$ be the corresponding sections over $\G_m^\ell \times \calC$. Then, there exists an integer $m=m(P,Q,\gamma)$ such that $\ord_\gamma D_{nP,Q} \leq m$ for every $n \geq 1$.   
	\end{lemma}
	\begin{proof}
		Since the statement is local we can assume, as in Lemma \ref{lemma:Betti_ab}, that the section $\sigma_{nP-Q}$ is given locally as
		\[
			\sigma_{nP-Q}: t \mapsto (h_1(t),\dots,h_\ell(t),t),
		\]
		where the $h_j$ are rational functions.
		Then, the order satisfies 
		\[
			\ord_\gamma D_{nP-Q} = \min_i \left\{ \dfrac{d^{(i)}h_j}{dt}(\gamma) \neq 0 \text{ for some	 } j=1\dots,\ell\right\}.
			\]
		 In particular, we can bound $\ord_\gamma D_{nP-Q}$ in terms of the minimal order of the derivative of $h_j$ that, for a fixed $j$, does not vanish at $\gamma$. This shows that we can reduce to the case in which $\ell=1$ so that we are considering the problem for $\G_m$. In this case, 
		\[
		h_1(t)=\frac{a^n(t)d(t)-b^n(t)c(t)}{b^n(t)c(t)},
		\] 
		where $a(t),b(t)$ and $c(t),d(t)$ are pairs of coprime polynomials with $P={a(t)}/{b(t)}$ and $Q={c(t)}/{d(t)}$ (recall that the operation of the group in this case is the multiplication). Therefore, it is enough to bound the multiplicities of the zeros of the numerator of $h_1$. We claim that we can conclude by applying a result of Ostafe \cite[Lemma 2.9]{Ostafe}. In order to see this, we notice that, by the coprimality assumption, the multiplicities of the common zeros of $a^n d$ and $b^n c$ are bounded independently of $n$. Hence we can factor them out and reduce to bound the multiplicities of $a'^{n}(t)d'(t)- b'^n(t)c'(t)$, where $a'd'$ and $b'c'$ have no common zeros; in particular, they satisfy the hypotheses of \cite[Lemma 2.9]{Ostafe} as wanted.
	\end{proof}

	Combining Lemma \ref{lemma:mult_ab} and Lemma \ref{lemma:mult_Gm} we can now prove Proposition \ref{prop:ord_boundQ}. 

	\begin{proof}[Proof of Proposition \ref{prop:ord_boundQ}]

	Given the split semiabelian scheme $\calG= \calA \times \G_m^{\ell} \rightarrow \calC$ $D_{nP,Q}$ we denote by $\pi_1$ and $\pi_2$ the projections from $\calG$ to $calA$ and $\G_m^{\ell}$ respectively. Given $D_{nP,Q}$, we denote by $D_{nP,Q}^\calA$ the divisor corresponding to the sections $\pi_1 \circ \sigma_P$ and $\pi_1 \circ \sigma_Q$ in the abelian scheme $\calA \to \calC$. Similarly, we denote by $D_{nP,Q}^{\G_m^\ell}$ the divisor corresponding to the sections $\pi_2 \circ \sigma_P$ and $\pi_2 \circ \sigma_Q$ in $\G_m^\ell \times \calC$.
	
	By the definition of $\ord_\gamma D_{nP,Q}$ and the fact that $\calG = \calA \times \G_m^\ell$, we see that 
	\begin{equation*}
		\ord_\gamma D_{nP,Q} \leq \min \{\ord_\gamma D_{nP,Q}^\calA , \ord_\gamma D_{nP,Q}^{\G_m^\ell }\}.
	\end{equation*} 

	Then the conclusion follows combining Lemma \ref{lemma:mult_ab} and Lemma \ref{lemma:mult_Gm}.
\end{proof}

\section{Proofs of Theorems \ref{thm:semiAb} and Theorem \ref{thm:PQ}}

We begin by proving Theorem \ref{thm:PQ}, which also implies part (i) of Theorem \ref{thm:semiAb}.

First, recall that $\gamma \in \supp D_{nP,Q}$ if and only if $\sigma_{nP}(\gamma)=\sigma_Q(\gamma)$. By Theorem \ref{thm:support}, we have that 
	\[ \calS := \bigcup_{n\ge 1} \supp D_{nP,Q}\]
	is a finite set. Moreover, by Proposition \ref{prop:ord_boundQ}, if $\ord_{\gamma} D_{nP,Q}\ge 1$, then the order is bounded independently of $n$. For every $\gamma \in \calS$ we put $m_{\gamma}:= \max_{n \ge 1} \ord_{\gamma} D_{nP,Q}$. Then,
	\[
	D:= \sum_{\gamma \in \calS} m_{\gamma} (\gamma),
	\]
	is the desired divisor.
	\medskip
	
	Let us now prove part (ii) of Theorem \ref{thm:semiAb}. Recall that $\ord_{\gamma} D_P\ge 1$ if and only if $\sigma_{P}(\gamma)=\calO(\gamma)$, and in this case by Lemma \ref{lemma:ord_bound} we have that $\ord_{\gamma} D_{nP}=\ord_{\gamma} D_P$ for every $n \neq 0$. It can however happen that $\ord_{\gamma} D_{nP}\ge 1$ but $\ord_{\gamma} D_{P}=0$; we are going to show that this can happen only for few choices of $n$.
	
	As before, by Theorem \ref{thm:support} we have that $\calS'= \cup_{n\ge 1} \supp D_{nP}$ is finite. Consider $\gamma \in \calS'\setminus  \supp D_{P}$ and let $n_{\gamma}$ be the smallest positive integer $n$ such that $\gamma \in  \supp D_{nP}$. Now, $\gamma \in \supp D_{nP}$ if and only if $n_{\gamma}$ divides $n$. Moreover, if $\gamma$ is not in the support of $D_P$ we have that $n_{\gamma}> 1$. This implies that for every positive integer $n$ not divisible by any of these finitely many $n_{\gamma}$ we have that $D_{nP}=D_P$, concluding the proof.

\begin{remark}
It is worth noticing that the statement (ii) of Theorem \ref{thm:semiAb} does not hold in general in the setting of Theorem \ref{thm:PQ}, i.e. when $\sigma_Q$ is not identically torsion. Indeed, it may happen that the two sections $\sigma_P$ and $\sigma_Q$ intersect at a point $\gamma \in \calC(\Qbar)$ and $\sigma_P(\gamma)$ is not torsion in $\calG_{\gamma}$. This implies that $\gamma \in \supp D_{P,Q}$ but $\gamma \not \in \supp D_{nP,Q}$ for every $n\geq 2$. 
\end{remark}



\bibliographystyle{amsalpha}	
\bibliography{gcd}

\end{document}